\documentclass{amsproc}

\usepackage{a4wide}
\usepackage{enumerate}

\newtheorem{Cor}{Corollary}

\newcommand{\tx}{\tilde{x}}
\renewcommand{\d}{{\mathrm d}}
\newcommand{\bt}{\begin{theorem}\ \ }
\newcommand{\et}{\end{theorem}}
\newcommand{\bp}{\begin{Prop}\ \ }
\newcommand{\ep}{\end{Prop}}
\newcommand{\bc}{\begin{Cor}\ \ }
\newcommand{\ec}{\end{Cor}}
\newcommand{\bl}{\begin{lemma}\ \ }
\newcommand{\el}{\end{lemma}}
\newcommand{\bd}{\begin{definition}\ \ }
\newcommand{\ed}{\end{definition}}
\newcommand{\pf}{\begin{proof}}
\newcommand{\epf}{\end{proof}}
\newcommand{\br}{\begin{remark}\ \ }
\newcommand{\er}{\end{remark}}
\newcommand{\brsn}{\begin{remarks*}\ \ }
\newcommand{\ersn}{\end{remarks*}}

\newcommand{\g}{\mathfrak{g}}
\newcommand{\h}{\mathfrak{h}}

\newcommand{\Ric}{\mathrm{Ric}}

\newcommand{\arr}{\begin{array}{rlll}}
\newcommand{\ea}{\end{array}}
\newcommand{\bea}{\begin{eqnarray}}
\newcommand{\eea}{\end{eqnarray}}
\newcommand{\bean}{\begin{eqnarray*}}
\newcommand{\eean}{\end{eqnarray*}}

\usepackage{
amsfonts, amssymb,
amsthm, amsgen
}
\newcommand{\cal}{\mathcal}   

\newcommand{\bcase}{\begin{case}}
\newcommand{\ecase}{\end{case}}

\newcommand{\bclaim}{\begin{claim}}
\newcommand{\eclaim}{\end{claim}}

\newcommand{\bstep}{\begin{step}}
\newcommand{\estep}{\end{step}}

\newcommand{\bhlem}{\begin{hlem}}
\newcommand{\ehlem}{\end{hlem}}

\newcommand{\bleer}{\begin{leer}}
\newcommand{\eleer}{\end{leer}}
\newcommand{\bde}{\begin{definition}}
\newcommand{\ede}{\end{definition}}

\newcommand{\bs}{\begin{prop}}
\newcommand{\es}{\end{prop}}
\newcommand{\btheo}{\begin{theorem}}
\newcommand{\etheo}{\end{theorem}}
\newcommand{\bfolg}{\begin{Cor}}
\newcommand{\efolg}{\end{Cor}}
\newcommand{\blem}{\begin{lem}}
\newcommand{\elem}{\end{lem}}
\newcommand{\bnote}{\begin{note}}
\newcommand{\enote}{\end{note}}
\newcommand{\bprf}{\begin{proof}}
\newcommand{\eprf}{\end{proof}}
\newcommand{\be}{\begin{eqnarray*}}
\newcommand{\ee}{\end{eqnarray*}}
\newcommand{\beqa}{\begin{eqnarray}}
\newcommand{\eeqa}{\end{eqnarray}}
\newcommand{\bi}{\begin{itemize}}
\newcommand{\ei}{\end{itemize}}
\newcommand{\bnum}{\begin{enumerate}}
\newcommand{\enum}{\end{enumerate}}

\newcommand{\beq}{\begin{equation}}
\newcommand{\eeq}{\end{equation}}

\newcommand{\vf}{\varphi}
\newcommand{\barr}{\[\begin{array}}
\newcommand{\earr}{\end{array}\]}
\newcommand{\bvec}{\left(\begin{array}{c}}
\newcommand{\evec}{\end{array}\right)}

\newcommand{\lah}{\mathfrak{h}}

\newcommand{\+}{\oplus}

\newcommand{\lason}{{\mathfrak{so}(n)}}

\newcommand{\p}{\partial}
\newcommand{\tp}{\tilde{\partial}}

\newcommand{\bbem}{\begin{remark}}
\newcommand{\ebem}{\end{remark}}
\newcommand{\bbez}{\begin{bez}}
\newcommand{\ebez}{\end{bez}}
\newcommand{\bbsp}{\begin{bsp}}
\newcommand{\ebsp}{\end{bsp}}

\newtheorem{theorem}{Theorem}
\newtheorem{prop}{Proposition}

\newtheorem{lem}{Lemma}

\theoremstyle{remark}
\newtheorem{remark}{Remark}

\newtheorem{ex}{Example}

\def\Ric{\mathop{{\rm Ric}}\nolimits}
\def\tRic{\mathop{\widetilde{\rm Ric}}\nolimits}

\def\Real{\mathbb{R}}
\def\g{\mathfrak{g}}
\def\h{\mathfrak{h}}

\def\so{\mathfrak{so}}
\def\spin{\mathfrak{spin}}
\def\simil{\mathfrak{sim}}

\def\su{\mathfrak{su}}
\def\u{\mathfrak{u}}
\def\sp{\mathfrak{sp}}

\def\zr{\ltimes}

\def\tr{\mathop\text{\rm tr}\nolimits}

\def\Ga{\Gamma}

\def\P{\mathcal{P}}

\begin{document}

 \bibliographystyle{abbrv}

 \title[On the local structure of Lorentzian Einstein  manifolds] {On the local structure of Lorentzian Einstein manifolds with parallel distribution of null lines}

\author{Anton S. Galaev}
 \address[Galaev]{Department of Mathematics and Statistics, Faculty of Science, Masaryk University in Brno, Kotl\'a\v rsk\'a~2, 611~37 Brno,
Czech Republic } \email{galaev@math.muni.cz}
\author{Thomas Leistner}
\address[Leistner]{School of Mathematical Sciences, The University of Adelaide, SA 5005,Australia}
 \email{thomas.leistner@adelaide.edu.au}

 \keywords{Einstein manifolds, Lorentzian manifolds, special holonomy.\\
 Version of \today}


\begin{abstract}
We study transformations of coordinates on a Lorentzian Einstein manifold with a parallel distribution of null
lines and show that the general Walker coordinates can be simplified.
In these coordinates, the full Lorentzian Einstein equation is reduced to  equations on a family of  Einstein  Riemannian  metrics.
\end{abstract}

\maketitle

\centerline{Dedicated to Dmitri Vladimirovich Alekseevsky on his
70th birthday}

\setcounter{tocdepth}{2}





\section{Introduction and statement of results}

Recently in \cite{G-P} G.W.~Gibbons and  C.N.~Pope considered the
Einstein equation on Lorentzian manifolds with holonomy algebras
contained in $\simil(n)$. A Lorentzian manifold $(M,g)$ of
dimension $n+2$ has holonomy algebra contained in $\simil(n)$ if
and only if it admits a parallel distribution $l$ of null lines
($l$ is a vector subbundle of rank one of the tangent bundle of
$M$ such that it holds $g(X,X)=0$ and $\nabla_YX$ is a section of
$l$ for all sections $X$ of $l$ and vector fields $Y$ on $M$, here
$\nabla$ is the Levi-Civita connection defined by $g$). Lorentzian
manifolds with this property have special Lorentzian holonomy  and
are of interest both in geometry (e.g.
\cite{B-M,Bazajkin,Boubel,Bryant,Schimming,Walker}) and
theoretical physics (e.g. \cite{BCH1,CGHP,CFH,FF00,Gibbons09}).
Any such manifold admits local coordinates $x^+,x^1,
...,x^n,x^{-}$, the so-called {\em Walker coordinates},  such that
the metric $g$ has the form
\begin{equation}\label{Walker} g=2\d x^+\d x^{-}+h+2A\d x^{-}+H(\d x^{-})^2,\end{equation}
where $h=h_{ij}(x^1,...,x^n,x^{-})\d x^i\d x^j$ is an
$x^{-}$-dependent family of Riemannian metrics,  $A=A_i(x^1, \ldots , x^n, x^-)\ \d x^i$ is an $x^-$-dependent family of one-forms, and $H$ is a local
function on $M$, \cite{Walker}. The vector field $\p_{+}:=\frac{\p}{\p x^+}$ defines the
parallel distribution of null lines.   We assume that the
indices $i,j,k,\ldots $ run from $1$ to $n$, and the indices $a,b,c, \ldots$   run in $+, 1, \ldots,  n, -$   and we use the
Einstein convention for sums. Furthermore, given coordinates $(x^+, x^1, \ldots , x^n, x^-)$ or
$(\tx^+, \tx^1, \ldots , \tx^n, \tx^-)$ we
 write $\p_a:=\frac{\p}{\p x^a}$ and $\tilde{\p}_a:= \frac{\p}{\p \tx^a}$.

The Einstein equation is the fundamental equation
of General Relativity. In the absence of
matter it has the form
\begin{equation}\label{Einstein} \Ric=\Lambda g,\end{equation} where $g$
is a Lorentzian metric on a manifold $M$, $\Ric$ is the Ricci
tensor of the metric $g$, i.e. $\Ric_{a b}=R^c_{\ a b c}$, where $R$ is the curvature tensor of the metric $g$,  and $\Lambda \in \Real$ is the
cosmological, or Einstein constant. If a metric $g$ of a smooth  manifold
$(M,g)$ satisfies this equation, then $(M,g)$ is called {\it an
Einstein manifold}. If moreover $\Lambda=0$,
 then it is called
{\it vacuum Einstein} or {\em Ricci-flat}. In dimension 4 examples of Einstein metrics are constructed in \cite{GT,KG1,KG2,Lewandowski,Petrov,St}.

We assume that $n\geq 2$, since for $n=0$ the problem is trivial
and for $n=1$ the metric \eqref{Walker} cannot be  non-flat and
Einstein \cite{G-P,GalEHol}.

In \cite{G-P} it is shown
that the Einstein equation  for a Lorentzian  metric of the form \eqref{Walker}
implies
\begin{equation}\label{GP4.7} H=\Lambda\cdot
(x^+)^2+x^+H_1+H_0,\end{equation} where $H_0$ and $H_1$ do not
depend on $x^+$. Furthermore, in \cite{G-P} it is proved that Equation
\eqref{Einstein} is equivalent to Equation \eqref{GP4.7}  and the following system of
equations
 \begin{align} \Delta H_0-\frac{1}{2}F^{ij}F_{ij}-2A^i\partial_{i}H_1-H_1\nabla^iA_i+2\Lambda A^iA_i-2\nabla^i\dot{A}_i&\nonumber \\ \label{GP4.8}
+\frac{1}{2}\dot h^{ij}\dot h_{ij}+h^{ij}\ddot h_{ij}+\frac{1}{2} h^{ij}\dot h_{ij}H_1& =0,\\ \label{GP4.9} \nabla^jF_{ij}+\partial_{i}H_1-2\Lambda A_i
+\nabla^j\dot h_{ij}-\p_{i}(h^{jk}\dot h_{jk})&=0,\\
\label{GP4.10} \Delta  H_1-2\Lambda \nabla^i A_i+\Lambda h^{ij}\dot h_{ij}&=0,\\ \label{GP4.11} \Ric_{ij}&=\Lambda
h_{ij},\end{align}
where
$\Delta H_0=h^{ij}(\p_{i}\p_{j}H_0-\Gamma^k_{ij}\p_k H_0)$ is the Laplace-Beltrami operator of the metrics $h(x^-)$ applied to $H_0$,
$F_{ij}=\p_{i}A_j-\p_{j}A_i$ are the components of the differential of the one-form $A(x^-)=A_i\d x^i$. A  dot denotes the derivative with respect to $x^{-}$ and
$\nabla^i\dot{A}_i=(\nabla_{\p_{i}}(\dot A)^{\#})^i$ is the divergence w.r.t. $h(x^-)$ of $\dot{A}$.

 Of course, the Walker coordinates are not
defined canonically and any other Walker coordinates
$\tilde x^{+},\tx^1, \ldots , \tx^n,\tilde x^{-}$ such that $\tilde{\p}_{+}=\p_{{+}}$ are given
by the following transformation (see \cite{Schimming} and Section \ref{Proofs})
$$ \tilde x^{+}=x^{+}+\varphi(x^1,...,x^n,x^{-}),\quad
\tilde x^{i}=\psi^{i}(x^1,...,x^n,x^{-}),\quad \tilde
x^{-}=x^{-}+c.$$ Now, the aim of the paper  is to simplify these
coordinates on Einstein manifolds and, in consequence, find easier
equivalences to the Einstein equation when written in the new
coordinates. That the coordinates can be simplified in special
situations was already shown in \cite{Schimming}:
\begin{prop}[Schimming \cite{Schimming}]\label{schimmingprop2}
 Let $(M,g)$ be a Lorentzian manifold with a parallel null vector field.
Then there exist local coordinates $\left(U,(x^+,x^1, \ldots , x^n,x^-)\right)$ such that the metric is given as
$$ g= 2\d x^+ \d x^- + h_{kl} \d x^k\d x^l$$
 with $h_{kl}$ smooth functions on $U$ with $\p_+ h_{kl}=0$.
 \end{prop}

Note that the condition for $(M,g)$ to admit a parallel null
vector field is stronger then the condition to admit a parallel
distribution of null lines. The first result of the present paper
generalises Proposition \ref{schimmingprop2} to manifolds with
only a parallel distribution of null lines:

\btheo\label{theo1} Let $(M,g)$ be a Lorentzian manifold with a
parallel distribution of null lines. Then there exist local
coordinates $\left(U,(x^+,x^1, \ldots , x^n, x^-)\right)$ such
that the metric is given as
$$ g=\left( 2\d x^+ + H \d x^-  \right)\d x^- + h_{kl} \d x^k\d x^l$$
 with $H$ and $h_{kl}$  smooth functions on $U$ with $\p_+ h_{kl}=0$.
 \etheo

With respect to coordinates as in Theorem \ref{theo1} the Einstein
equations (\ref{GP4.8}--\ref{GP4.11}) become much easier:

\begin{align}\label{GP4.8CC}
\Delta H_0+\frac{1}{2}\dot h^{ij}\dot h_{ij}+h^{ij}\ddot
h_{ij}+\frac{1}{2} h^{ij}\dot h_{ij}H_1& =0,\\ \label{GP4.9CC}
\partial_{i}H_1+\nabla^j\dot h_{ij}-\p_{i}(h^{jk}\dot h_{jk})&=0,
\\
\label{GP4.10CC} \Delta  H_1+\Lambda h^{ij}\dot h_{ij}&=0,\\
\label{GP4.11CC} \Ric_{ij}&=\Lambda h_{ij}.\end{align}

Then we assume that the manifold is Einstein, and, based on
Equation \eqref{GP4.7}, we prove the following:

\btheo\label{theo2} Let $(M,g)$ be a Lorentzian manifold with a
parallel distribution of null lines and assume that $M$ is
Einstein with Einstein constant $\Lambda$. Then there exist local
coordinates $\left(x^+,x^1, \ldots , x^n, x^-\right)$ such that
the metric is given as
$$ g=\left( 2\d x^+ + (\Lambda (x^+)^2 + x^+ H_1 ) \d x^-  \right)\d x^- +  h_{kl} \d x^k\d x^l$$
 with $H_1$ and $h_{kl}$  smooth functions on $U$ with $\p_+ h_{kl} =\p_+H_1=0$ and satisfying the equations
 \begin{align}\label{GP4.8C}
\frac{1}{2}\dot h^{ij}\dot h_{ij}+h^{ij}\ddot h_{ij}+\frac{1}{2} h^{ij}\dot h_{ij}H_1& =0,\\ \label{GP4.9C} \partial_{i}H_1+\nabla^j\dot h_{ij}-\p_{i}(h^{jk}\dot h_{jk})&=0,
\\
\label{GP4.10C} \Delta  H_1+\Lambda h^{ij}\dot h_{ij}&=0,\\ \label{GP4.11C} \Ric_{ij}&=\Lambda
h_{ij}.\end{align}
Conversely, any such metric is Einstein with Einstein constant $\Lambda$.
 \etheo


Note that if $(M,g)$ admits a parallel null vector field, then the
Walker coordinates in \eqref{Walker} satisfy  $\p_{+} H=0$ and we
get Proposition \ref{schimmingprop2}  from Theorem \ref{theo1}
(see Remark \ref{rem1} below). If, in addition,  such a metric is
Einstein, then $\Lambda=0$, i.e. this metric is Ricci-flat and the
equations (\ref{GP4.8}--\ref{GP4.11}) take the following more
simplified form
\begin{align}\label{GP4.8B} \frac{1}{2}\dot h^{ij}\dot h_{ij}+h^{ij}\ddot h_{ij}& =0,\\ \label{GP4.9B} \nabla^j\dot h_{ij}-\p_{i}(h^{jk}\dot h_{jk})&=0,\\
 \label{GP4.11B} \Ric_{ij}&=0.\end{align}

Finally, as the main result of the paper we show that in the the
case $\Lambda\neq 0$ we can do better.

\begin{theorem} \label{Main Theorem} Let $(M,g)$ be a Lorentzian  manifold of
dimension $n+2$
 admitting a parallel distribution of null lines. If $(M,g)$ is Einstein with the non-zero
 cosmological constant $\Lambda$ then  there exist local coordinates
$\left(x^+,x^1, \ldots , x^n, x^{-}\right)$ such that the metric $g$ has the form
$$ g=2\d x^+\d x^{-}+h_{kl}\d x^k\d x^l+(\Lambda
(x^+)^2+H_0)(\d x^{-})^2$$
with $\p_+h_{kl}=\p_+ H_0=0$,  $h_{kl}$ defines
 an $x^-$-dependent family of  Riemannian  Einstein metrics with the cosmological
constant $\Lambda$,  satisfying the equations
\begin{align} \label{GP4.8A} \Delta
H_0+\frac{1}{2}h^{ij}\ddot{h}_{ij} &=0,\\
\label{GP4.9A}\nabla^j\dot h_{ij}&=0, \\
\label{GP4.10A}h^{ij}\dot{h}_{ij} &=0, \\
\label{GP4.11A}\Ric_{ij}&=\Lambda h_{ij},\end{align}
where
$\dot{h}_{ij}=\partial_{-}{h_{ij}}$.
Conversely,
any such metric is Einstein. \end{theorem}

Remark that in \cite{G-P} it is shown that Equation \eqref{GP4.10}
follows from \eqref{GP4.9} and \eqref{GP4.11}, i.e. it may be
omitted from the Einstein equation. By the same reason Equations
\eqref{GP4.10CC} \eqref{GP4.10C} and \eqref{GP4.10A} may be
omitted. On the other hand, these equations can be used as the
corollaries of the Einstein equation.

Thus, we reduce the Einstein equation with $\Lambda\neq 0$ on a
Lorentzian manifold with holonomy algebra contained in $\simil(n)$
to the study of families of Einstein  Riemannian metrics
satisfying Equation \eqref{GP4.9A}. If $\Lambda=0$ and $\p_+ H\neq
0$, i.e. $H_1\neq 0$, then consider the coordinates as in Theorem
\ref{theo2}. Equation \eqref{GP4.10C} shows that $H_1$ is a family
of harmonic functions on the family of the Riemannian manifolds
with metrics $h(x^-)$. Fixing any such $H_1$ we get Equations
\eqref{GP4.8C} and \eqref{GP4.9C} on the family of Ricci-flat
Riemannian  metrics $h(x^-)$. Finally, if $(M,g)$ is Einstein and
it admits a parallel null vector field, then it is Ricci flat and
this is equivalent to Equations   \eqref{GP4.8B} and
\eqref{GP4.9B} on the family of Ricci-flat Riemannian  metrics
$h(x^-)$. In Section \ref{secConseq}  we consider the holonomy
algebra of $(M,g)$ and the de~Rham decomposition for the family of
Riemannian metrics $h(x^-)$.

Note that to find the required  transformation of the coordinates in Theorem \ref{Main Theorem}
we need to solve a system of ODE's,
while in  \cite{Schimming} several  PDE's need to be solved.

Examples of Einstein metrics of the form \eqref{Walker} with $h$
independent of $x^-$ and each possible holonomy algebra are
constructed in \cite{GalEHol}. It is interesting to construct
examples of Einstein manifolds satisfying some global properties,
e.g. global hyperbolic, as in \cite{B-M} or \cite{Bazajkin}.

In Section \ref{secEx} we consider examples in dimension 4.

\section{Consequences}\label{secConseq}

Let us consider some consequences of the above theorems. Let $(M,g)$ be a Lorentzian manifold with a parallel distribution of null lines.
Without loss of generality we may assume that $(M,g)$ is locally indecomposable, i.e. locally
it is not a product of a Lorentzian and of a Riemannian manifold. The holonomy of such manifolds are contained in
$\simil(n)=(\Real\oplus\so(n))\zr \Real^n$. In \cite{Leistner} it was shown that the projection $\lah$  of the holonomy algebra of $(M,g)$ onto $\lason$ has to be a Riemannain holonomy algebra.
Now,
recall that for each Riemannian holonomy algebra $\h\subset\so(n)$ there exists a  decomposition
\begin{equation}\label{LM0A}\Real^{n}=\Real^{n_0}\oplus\Real^{n_1}\oplus\cdots\oplus\Real^{n_r}\end{equation} and the
corresponding decomposition into the direct sum of ideals
\begin{equation}\label{LM0B}\h=\{0\}\oplus\h_1\oplus\cdots\oplus\h_r\end{equation}
such that each  $\h_\alpha\subset\so(n_\alpha)$ is an irreducible Riemannian  holonomy algebra,
in particular it coincides with one of the following  subalgebras of $\so(n_\alpha)$:  $\so(n_\alpha)$, $\u(\frac{n_\alpha}{2})$, $\su(\frac{n_\alpha}{2})$,
$\sp(\frac{n_\alpha}{4})\oplus\sp(1)$, $\sp(\frac{n_\alpha}{4})$, $G_2\subset\so(7)$, $\spin_7\subset\so(8)$
or it is an irreducible symmetric Berger algebra (i.e. it is the
holonomy algebra of a symmetric Riemannian manifold and it is
different from $\so(n_\alpha)$, $\u(\frac{n_\alpha}{2})$, $\sp(\frac{n_\alpha}{4})\oplus\sp(1)$). Recall that if the holonomy algebra of a  Riemannian manifold is a symmetric
Berger algebra, then the manifold is locally symmetric.

In \cite{ESI,GalEHol} it is proven that if $(M,g)$ is Einstein with $\Lambda\neq 0$, then the holonomy algebra of $(M,g)$ has the form
$\g=(\Real\oplus\h)\zr\Real^n$, moreover, each subalgebra $\h_\alpha\subset\so(n_\alpha)$ from the decomposition \eqref{LM0B} coincides with one of the algebras
$\so(n_\alpha)$, $\u(\frac{n_\alpha}{2})$, $\sp(\frac{n_\alpha}{4})\oplus\sp(1)$ or with a symmetric Berger algebra,  and in the  decomposition \eqref{LM0A} it holds $n_0=0$.
Next, if $\Lambda=0$, then one of the following holds:

\begin{itemize}\item[(A)] $\g=(\Real\oplus\h)\zr\Real^n$ and at least one of the subalgebras $\h_\alpha\subset\so(n_\alpha)$
from the decomposition \eqref{LM0B} coincides with one of the algebras
$\so(n_\alpha)$, $\u(\frac{n_\alpha}{2})$, $\sp(\frac{n_\alpha}{4})\oplus\sp(1)$ or with a symmetric Berger algebra.
\item[(B)] $\g=\h\zr\Real^n$  and each subalgebra $\h_\alpha\subset\so(n_\alpha)$ from the decomposition \eqref{LM0B}
coincides with one of the algebras $\su(\frac{n_\alpha}{2})$,  $\sp(\frac{n_\alpha}{4})$, $G_2\subset\so(7)$, $\spin_7\subset\so(8)$.
\end{itemize}

In \cite{Boubel} it is proved that there exist Walker coordinates
$x^+,x_0^1,\ldots ,x_0^{n_0},\ldots ,x_r^1,...,x_r^{n_r},x^-$ that
are adapted to the decomposition \eqref{LM0B}. This means that
$h=h_0+h_1+\cdots+h_r,$ $h_0=\sum_{i=1}^{n_0}(\d x_0^i)^2$ and
$A=\sum_{\alpha=1}^r\sum_{k=1}^{n_\alpha} A^\alpha_k \d
x^k_\alpha$ and for each $1\leq \alpha\leq r$ it holds
$h_\alpha=\sum_{i,j=1}^{n_\alpha}h_{\alpha ij} \d x_{\alpha}^i\d
x_{\alpha}^j$ with $\frac{\p}{\p {x^k_\beta}}h_{\alpha
ij}=\frac{\p}{\p {x^k_\beta}}A^\alpha_i=0$ for all $1\leq i,j\leq
n_\alpha$ if $\beta\neq\alpha$. We will show that the
transformations can be chosen in such a way that the new
coordinates  are adapted in this sense.

\begin{prop}\label{prop2}
Let $(M,g)$ be a Lorentzian manifold with a parallel distribution
of null lines and let $\lah$ be the projection of its holonomy
algebra onto $\lason$ decomposed as in \eqref{LM0B}. \bnum
\item
Then the coordinates found in Theorem \ref{theo1}  can be chosen to be adapted to this decomposition.
\item If $(M,g)$ is Einstein with $\Lambda\neq 0$, then there exist coordinates adapted to this decomposition with the properties as in
Theorem \ref{Main Theorem} and with $n_0=0$.
\enum
\end{prop}
We will  prove this proposition in the next section. It shows that
the Einstein conditions written  as in the formulae after Theorem
\ref{theo1} and in Theorem \ref{Main Theorem} can, in addition,
be formulated in adapted coordinates.

 Now we discuss to which extend the Einstein equations in the theorems have to be satisfied for each of the $h_\alpha$'s separately when written in the coordinates of Proposition \ref{prop2}.
First, let $\Lambda\neq 0$ and consider (\ref{GP4.8A} -- \ref{GP4.11A}). It is obvious that each $h_\alpha$ satisfies
\eqref{GP4.9A} and \eqref{GP4.11A}. Using the first variation formula for the Ricci tensor (see e.g. \cite[Theorem 1.174]{Besse}),  in \cite{G-P} it is shown that \eqref{GP4.10} follows from \eqref{GP4.9} and
\eqref{GP4.11}
 by taking the divergence of \eqref{GP4.9}. Hence, using the divergence with respect to the metric $h_\alpha$,
 \eqref{GP4.9A} and \eqref{GP4.11A} imply that  each $h_\alpha$ satisfies also \eqref{GP4.10A}.
This means that one has to solve   (\ref{GP4.9A} -- \ref{GP4.11A}) separately for each $h_\alpha$
and then find $H_0$ from \eqref{GP4.8A}.

Similarly, if $\Lambda=0$, consider (\ref{GP4.8CC} --
\ref{GP4.11CC}). Obviously, each $h_\alpha$ has to be Ricci-flat.
Applying the divergence with respect to $h_\alpha$ to
\eqref{GP4.9C} we get that $\Delta_\alpha H_1=0$. This together
with \eqref{GP4.9C} shows that $H_1=\sum_\alpha H_{1\alpha}$,
where each $H_{1\alpha}$ depends only on $x^i_\alpha$  and it is
harmonic with respect to $h_\alpha$. Now each $h_\alpha$ satisfies
\eqref{GP4.9CC} with $H_1$ replaced by $H_{1\alpha}$.

Next we study the possible  summands in the decomposition
\eqref{LM0B} under the assumption that the manifold $(M,g)$ is
Einstein with $\Lambda\not=0$. First we claim that if
$\h_\alpha\subset\so(n_\alpha)$ is a symmetric Berger algebra,
then  each metric in the family $h_\alpha(x^{-})$ is locally
symmetric and its holonomy algebra coincides with $\h_\alpha$.
Indeed, the  holonomy algebra $\h_\alpha(x^{-})$ of  each metric
in the family $h_\alpha(x^{-})$ is contained in $\h_\alpha$ and it
is non-trivial due to \eqref{GP4.11A}. Since
$\h_\alpha\subset\so(n_\alpha)$ is a symmetric Berger algebra its
space of curvature tensors ${\cal R}(\h_\alpha)$ is
one-dimensional. This shows that $\h_\alpha(x^{-})=\h_\alpha$. If
$\h_\alpha=\u(\frac{n_\alpha}{2})$ (resp.,
$\h_\alpha=\sp(\frac{n_\alpha}{4})\oplus\sp(1)$), then each metric
in the family $h_\alpha(x^{-})$ is K\"ahler-Einstein (resp.,
quaternionic-K\"ahler). For some values of $x^{-}$ the metric
$h_\alpha(x^{-})$ can be decomposable, but it does not contain a
flat factor. If $\h_\alpha=\so(n_\alpha)$, then we get a general
family of Einstein metrics. For some values of $x^{-}$ the metric
$h_\alpha(x^{-})$ can be decomposable, but it does not contain a
flat factor.
\begin{prop}\label{prop1} Under the current assumptions, if $\h_\alpha\subset\so(n_\alpha)$ is a symmetric Berger algebra,
 then $h_\alpha$ satisfies the equation
 \begin{equation}\label{4.9strong} \nabla_{i}(h_\alpha^{kt}\dot h_{\alpha tj})-2\dot \Gamma^k_{ij}=0,\quad 1\leq i,j,k\leq n_\alpha,
\end{equation}
where $\Gamma^k_{ij}$ is the family of  the Christoffel symbols
for the family of the Riemannian metrics  $h(x^{-})$. \end{prop}

 Note that the Equation \eqref{4.9strong} is stronger then the
Equation \eqref{GP4.9A}, since the last equation is obtained from
the first one by taking the trace. This proposition will be proved
below.

Finally, suppose that $\Lambda=0$ and the holonomy algebra $\g$ of
$(M,g)$ is as in the case (A) above. Suppose that $\h_\alpha$ is
one of $\so(n_\alpha)$, $\u(\frac{n_\alpha}{2})$,
$\sp(\frac{n_\alpha}{4})\oplus\sp(1)$ or it is  a symmetric Berger
algebra. Equation \eqref{GP4.11C}  shows that in the first three
cases each metric in the family $h_\alpha$ is Ricci-flat,
consequently, its holonomy algebra is contained, respectively, in
$\so(n_\alpha)$, $\su(\frac{n_\alpha}{2})$,
$\sp(\frac{n_\alpha}{4})$. If $\h_\alpha$ is a symmetric Berger
algebra, then by the same reasons each metric in the family
$h_\alpha (x^-)$ is flat. Otherwise $\h_\alpha$ is either trivial
or it is  one of  $\su(\frac{n_\alpha}{2})$,
$\sp(\frac{n_\alpha}{4})$, $G_2\subset\so(7)$,
$\spin_7\subset\so(8)$. Each  metric in the family $h_\alpha$ is
Ricci-flat and it has holonomy algebra contained in $\h_\alpha$.

To sum up the consequences
we remark that
 the problem of finding Einstein Lorentzian metrics with
$\Lambda\neq 0$ is reduced first to the problem of finding
families of Einstein Riemannian metrics satisfying Equation
\eqref{GP4.9A} (or \eqref{4.9strong} for the symmetric case) and
then to Poisson equation \eqref{GP4.8A} on the function $H_0$.
This is related to the module spaces of Einstein metrics
\cite{Besse}. For example, for most of symmetric Berger algebras
$\h_\alpha\subset\so(n_\alpha)$ it holds that $h_\alpha$ is an
isolated metric, i.e. it  is independent of $x^{-}$
\cite{Koi2,Koi4}. Hence, since it is symmetric, it is uniquely
defined by $\Lambda$. Similarly, if $\Lambda=0$ and $\p_+ H\neq
0$, i.e. $H_1\neq 0$, then consider the coordinates as in Theorem
\ref{theo1}. Equation \eqref{GP4.10CC} shows that $H_1$ is a
family of harmonic functions on the family of the Riemannian
manifolds with metrics $h(x^-)$. Fixing any such $H_1$ we get
Equation \eqref{GP4.9CC} on the family of Ricci-flat Riemannian
metrics $h(x^-)$ and then  Poisson equation \eqref{GP4.8CC} on the
function $H_0$. Finally, if $(M,g)$ is Einstein and it admits a
parallel null vector field, then it is Ricci-flat and this is
equivalent to the equations \eqref{GP4.8B} and \eqref{GP4.9B} on
the family of Ricci-flat Riemannian  metrics $h(x^-)$.

\section{Proofs}
\label{Proofs}


\subsection*{Coordinate transformations}
In order to simplify the Walker coordinates, first we
have to describe the most general coordinate transformation leaving the form \eqref{Walker} invariant.
This was already done in \cite{Schimming} in the case of a parallel null vector field.

\begin{prop}\label{schimmingprop1}
The most general coordinate transformation  with
$\tilde{\p}_+=\p_+$ that  preserves the form \eqref{Walker} is
given by
\begin{equation} \tilde x^{+}=x^{+}+\varphi(x^1,...,x^n,x^{-}),\quad
\tilde x^{i}=\psi^{i}(x^1,...,x^n,x^{-}),\quad \tilde
x^{-}=x^{-}+c. \label{gauge}
\end{equation}
If  the metric and its inverse is written as
\begin{equation}\label{metric}
 g=\left(\begin{array}{ccc}
 0&0& 1 \\ 0 & h & A\\
 1& A^t & H \end{array}\right)
 \ \text{ and }\
 g^{-1}=\left(\begin{array}{ccc}
 F&B^t& 1 \\ B & h^{-1} & 0\\
 1& 0 & 0 \end{array}\right),
 \end{equation}
 with $B=-h^{-1}A$ and $F+H+A^tB=0$, then in the
 new coordinates  it holds
\begin{eqnarray}
\label{gtrans} \tilde{h}^{ij}&=&\p_k \psi^i h^{kl} \p_l\psi^{j} \\
\label{Btrans} \tilde{B}^i
 &=& \p_-\psi^i + B^k \p_k\psi^i +  h^{kl} \p_k\vf \p_l \psi^i \\
\label{Ftrans}  \tilde{F}&=&
F+ \p_- \vf +  B^k  \p_k \vf
 +
 h^{kl} \p_k \vf \p_l \vf.
\end{eqnarray}
\end{prop}
\bprf Since $\tilde{\p}_+=\p_+$,  the transformation formula for
the canonical basis implies
$$ \p_+\tx^-  = 0,\ \p_+ \tx^k  = 0, \text{ and }\ \p_+\tx^+=1.$$
Furthermore we get
$$0=g(\p_+,\p_i)=
\p_+\tx^+ \p_i\tx^- g(\tp_+,\tp_-) + \p_+\tx^+ \p_i \tx^k g(\tp_+,\tp_k).
$$
As we require $g(\tp_+,\tp_k)=0$ this implies $\p_i\tx^-=0$.
Finally we have to check
$$ 1=g(\p_+,\p_-)=
\p_+\tx^+\p_-\tx^-,$$
which implies
 $\p_-\tx^-=1$.
This shows that the most general transformation is of the form \eqref{gauge}.

 In order to write down the inverse metric coefficients in the new  coordinates first we see that
 in the coordinates $\eqref{Walker}$   the metric and its inverse are given as in \eqref{metric}.
 The transformation formula for the inverse metric coefficients $g^{ab}$ is given by
 $$
 \p_c \tx^a g^{cd} \p_d\tx^{b} = \tilde{g}^{ab},$$
 where $a$ and $b $ run over $+, 1, \ldots , n, -$.
 This implies that
 \be \tilde{B}^i\ =\ \tilde{g}^{+i}
 &=&
 \p_-\tx^i +  B^k \p_k\tx^i + h^{kl}\p_k\tx^+ \p_l \tx^i,
  \ee
  which is Equation \eqref{Btrans}.
Furthermore we get
 \be
 \tilde{F}\ =\ \tilde{g}^{++}&=&
 F+ \p_+\tx^+ \p_-\tx^+ +
 B^k \p_+\tx^+ \p_k \tx^+
 +
 h^{kl}\p_k \tx^+ \p_l \tx^+,
\ee
which is Equation \eqref{Ftrans}.
In the same way the equations for $\tilde{h}^{ij}$.
\eprf
\subsection*{Proof of Theorem \ref{theo1}}
 Setting $\tilde{B}^i$ to zero for each $i=1, \ldots , n$ in the transformation formula above we obtain  a linear PDE for  the function $\psi$
 \begin{equation}\label{ck1}
 \p_-\psi =-\left( B^k +
  h^{kl}\p_l\vf\right) \p_k \psi ,
 \end{equation}
 and we have to find $n$ linear independent solutions $\psi^1, \ldots , \psi^n$.
This problem can be solved for the following reasons: Fix the function $\vf=\vf (x^1, \ldots , x^n,x^-)$, e.g. $\vf\equiv 0$, and consider the characteristic  vector field of \eqref{ck1}
$$X:= \p_- +  \left( B^k +  h^{kl}\p_l\vf\right) \p_k.$$
 Obviously, Equation \eqref{ck1} is equivalent to the equation
 \begin{equation}\label{kwneu}
 X(\psi)=d\psi(X)=0.
 \end{equation}
We have $\left[ \p_+,X  \right] =0$. Hence, we find coordinates $(y^+ , y^1, \ldots , y^n, y^- ) $ such that
$$ \tfrac{\p}{\p y^+}= \p_+\  \text{ and }\  \tfrac{\p}{\p y^-} =X.$$
Now, any function $\psi=\psi(y^1, \ldots , y^n)$ satisfies Equation \eqref{kwneu}. Note that $\p_+ y^-=\p_+y^i=0$ and therefore also
$\p_+\psi=0$. Taking $n$ linear independent solutions gives us the required solutions $\psi^i$ of Equation \eqref{ck1} to build the new coordinate system.
 %
\hfill $\Box$

 \bbem \label{rem1}
In order to obtain Schimming's result of Proposition \ref{schimmingprop2}
one has to set $\tilde{H}$ to zero obtaining the additional equation
\begin{equation}\label{ck2}
\p_- \vf=- F-  B^k  \p_k \vf
 -
 h^{kl}\p_k \vf \p_l \vf
\end{equation} together with the linear Equation \eqref{ck1}.
Although Equation \eqref{ck2} cannot be written in the form $X(\vf)=0$, it can be solved using characteristics (see below).
\ebem
\bbem
Note that Schimming's result
 cannot be true only with the assumption of a parallel distribution of null lines: Since in this case $H$ and thus $F$ may depend on $x^+$ but $\vf$ does not, Equation \eqref{ck2} cannot be solved. In other words, the $x^+$-dependence of $H$ in general cannot be changed by these coordinate transformations.
But in case of Einstein metrics with arbitrary Einstein constant \nolinebreak $\Lambda$, Theorem \ref{theo2} shows that one can get rid of the part of $H$ that does not depend on $x^+$.
 \ebem
%
%

\subsection*{Proof of Theorem \ref{theo2}}
We fix coordinates $(x^+,x^1, \ldots , x^n,x^-)$ as in Theorem \ref{theo1} with $A_i=0$. Since $(M,g)$ is Einstein it holds that
$$ H= \Lambda (x^+)^2 +x^+ H_1+H_0,$$
where $\p_+ H_1= \p_+ H_0=0$.
 Now we try to find an appropriate coordinate transformation consisting of functions  $\vf$ and $\psi^i$ as in Proposition \ref{schimmingprop1}.
 First we consider the equation
 \begin{equation}\label{phiequ} \p_-\vf \ =\ H_0 -H_1\vf+ \Lambda \vf^2- h^{kl}\p_k\vf \p_l \vf .\end{equation}
 This equation can be solved by the method of characteristics (for details see for example \cite[Chapter 10, Section 1]{spivak5}). Since the $x^-$ derivative of $\vf$ is isolated, a characteristic is given by $(x^1, \ldots , x^n)\mapsto (x^1, \ldots x^n,0)$ and the parameter of the characteristic curves can be chosen to be $x^-$.
  Let  $\vf$ be a smooth solution of this equation.
With respect to this $\vf$ we consider the equation
\begin{equation}\label{psiequ}
 \p_- \psi \ =\ - h^{kl} \p_k\vf \p_l\psi .
 \end{equation}
 As in Theorem \ref{theo1}, we find $n$ linear
 independent solutions $\psi^1, \ldots , \psi^n$ to this equation.  Hence, in the new
 coordinates given as in \eqref{gauge} we still have  $\tilde{B}^k=0$.
Now, since $(M,g)$ is Einstein, it is
$$ \tilde{H}\ =\    \Lambda (\tx^+)^2 +\tx^+ \tilde{H}_1+\tilde{H}_0
\ =\ \Lambda (x^+)^2 + (2\Lambda \vf + \tilde{H}_1)x^+
+\tilde{H}_1\vf +\Lambda \vf^2+\tilde{H}_0 .$$ On the other hand,
from the transformation formula and $\tilde{B}^k=0$ we have \be
\tilde{H}\ =\ -\tilde{F} &=& - F - \p_- \vf
 -
h^{kl}\p_k \vf \p_l \vf
\\
&=&
\left(\Lambda (x^+)^2 + x^+ H_1+H_0\right) - \p_- \vf
 -
 h^{kl}\p_k \vf \p_l \vf. \ee Comparing these two equations
and differentiating w.r.t. $\p_+$ shows that $ (2\Lambda \vf +
\tilde{H}_1)\ = \ H_1$ and furthermore
$$ \Lambda \vf^2+\tilde{H}_0 +\tilde{H}_1\vf
\ =\
H_0 - \p_- \vf
 -
 h^{kl}\p_k \vf \p_l \vf.
$$
Hence, putting this together we get
$$
\tilde{H}_0\ =\
H_0 - \p_- \vf
 -
 h^{kl}\p_k \vf \p_l \vf
+\Lambda \vf^2  -H_1\vf.
$$
But since $\vf$ satisfies Equation \eqref{phiequ}, we obtain $\tilde{H}_0=0$ in the new coordinates.
\hfill $\Box$

\subsection*{Curvature tensors} For the proof of
Theorem \ref{Main Theorem} we need some algebraic preliminaries.
The tangent space to $M$ at any point $m\in M$ can be identified
with the Minkowski space $\Real^{1,n+1}$. Denote by $g$ the metric
on it.  Let $\Real p$ be the null line corresponding to the
parallel distribution. Let ${\cal R}(\simil(n))$ be the space of
algebraic curvature tensors of type $\simil(n)$, i.e. the space of
linear maps from $\Lambda^2\Real^{1,n+1}$ to $\simil(n)$
satisfying the first Bianchi identity.   The curvature tensor
$R=R_m$ at the point $m$  belongs to the space ${\cal
R}(\simil(n))$. The space ${\cal R}(\simil(n))$ is found in
\cite{Gal1,GalP}. We will review this result now.    Fix
a null vector $q\in\Real^{1,n+1}$ such that $g(p,q)=1$. Let
$E\subset \Real^{1,n+1}$ be the orthogonal complement to $\Real
p\oplus\Real q$, then $E$ is an Euclidean space. We get the
decomposition
\begin{equation}\label{decompT}\Real^{1,n+1}=\Real p\oplus E\oplus\Real q.\end{equation}
We will often write $\Real^n$ instead of  $E$.  Fixing a basis $X_1,...,X_n$ in $\Real^n$, we get that
\begin{equation}\label{matsim}\simil(n)=\left\{\left. \left (\begin{array}{ccc} a
&(GX)^t & 0\\ 0 & A &-X \\ 0 & 0 & -a \\
\end{array}\right)\right|\, a\in \Real,\,A \in \so(n),\, X\in \Real^n \right\}, \end{equation} where $G$ is the Gram matrix of the metric $g|_{\Real^n}$ with respect
to the basis $X_1,...,X_n$.
The above matrix can be identified with the triple $(a,A,X)$. We obtain the decomposition
$$\simil(n)=(\Real\oplus\so(n))\zr \Real^n.$$
For a subalgebra $\h\subset\so(n)$  consider the space
$$\P(\h)=\{P\in (\Real^n)^*\otimes
\h|g(P(x)y,z)+g(P(y)z,x)+g(P(z)x,y)=0\text{ for all } x,y,z\in
\Real^n\}.$$ Define the map $\tRic:\P(\h)\to\Real^n$,
$\tRic(P)=P^j_{ik}g^{ik}X_j$. It does not depend on the choice of
the basis $X_1,...,X_n$. The tensor $R\in{\cal R}(\simil(n))$ is
uniquely given by elements $\lambda\in\Real,v\in E,R_0\in{\cal
R}(\so(n)),P\in\P(\so(n)),T\in\odot^2E$ in the following way.
\begin{align*}
R(p,q)=&(\lambda,0,v),\qquad R(x,y)=(0,R_0(x,y),P(y)x-P(x)y),\\
R(x,q)=&(g(v,x),P(x),T(x)),\qquad R(p,x)=0 \end{align*}
 for all $x,y\in\Real^n$.
We write $R=R(\lambda,v,R_0,P,T)$. The Ricci tensor $\Ric(R)$ of
$R$ is given by  $\Ric(R)(X,Y)=\tr(Z\mapsto R(X,Z)Y)$  and it
satisfies
\begin{align}\label{Ric1} \Ric(p,q)=&-\lambda,\quad \Ric(x,y)=\Ric(R_0)(x,y),\\
\label{Ric2} \Ric(x,q)=&g(x,\tRic(P)-v),\quad \Ric(q,q)=\tr T. \end{align}

Let us take some other null vector $q'$ with $g(p,q')=1$. There
exists a unique vector $w\in E$ such that $q'=-\frac{1}{2}g(w,w)
p+w+q$.
 The corresponding
$E'$ has the form $E'=\{-g(x,w)p+x|x\in E\}$. We will consider the
map $x\in E\mapsto x'=-g(x,w)p+x\in E'$. Using this, we  obtain
that $R=R(\tilde\lambda,\tilde v,\tilde R_0,\tilde P,\tilde T)$. For example, it holds
$$\tilde \lambda=\lambda,\quad \tilde v=(v-\lambda w)',\quad
\tilde P(x')=(P(x)-R_0(x,w))',\quad \tilde R_0(x',y')z'=(R_0(x,y)z)'.$$ This shows that
using the change of $q$ we may get rid of $v$ or some times of
$P$. (For example, if $\h$ is a symmetric Berger algebra, i.e.
$\dim{\cal R}(\h)=1$,  and $R_0\neq 0$, then there exists $w\in E$
such that
$P(x)-R_0(x,w)=0$ for all $x$ \cite{GalP}, i.e. $\tilde P=0$.)  

\subsection*{Proof of Theorem \ref{Main Theorem}}
Consider the general Walker metric \eqref{Walker}.
Suppose that it is Einstein with $\Lambda\neq 0$. Then
$H=\Lambda(x^+)^2+x^+H_1+H_0$, where $H_0$ and $H_1$ are
independent of $x^+$ \cite{G-P}. Consider the vector fields
$$p=\p_{+},\quad X_i=\p_{i}-A_i\p_{+},\quad q=\p_{-}-\frac{1}{2}H\p_{+}.$$
Let $E\subset TM$ be the distribution generated by the vector
fields $X_i$. At each point $m$ we get $$T_mM=\Real p_m\oplus
E_m\oplus\Real q_m.$$ Then the curvature tensor $R$ is given by the
elements $\lambda,v,R_0,P,T$ as above but depending on the point.
Since the manifold is Einstein, we get $\lambda=-\Lambda$.

\begin{prop} For any $W\in \Ga(E)$ such that
$\nabla_{\p_{+}}W=0$ there exist new Walker coordinates $\tilde x^{a}$ such that
the corresponding vector field $q'$ has the form  $q'=-\frac{1}{2}g(W,W)p+W+q$.
\end{prop}

\bprf  Let us write $W=W^iX_i$.  Since
$\nabla_{\p_{+}}W=0$, we get that $\p_{+} W^i=0$.
We will find the inverse transformation $$x^+=\tilde x^{+},\quad
x^i=x^i(\tilde x^{1},...,\tilde x^n,\tilde x^{-}),\quad  x^{-}=\tilde x^{-}.$$ It holds
$$\tilde\p_{+}=\p_{+}, \quad
\tilde\p_{i}=\frac{\p{x^j}}{\p{\tilde x^{i}}}\p_{j},\quad
\tilde\p_{-}=\frac{\p{x^i}}{\p{\tilde x^{-}}}\p_{i}+\p_{{-}}.$$
For the new Walker metric we have
$$H'=g(\tilde\p_{-},\tilde\p_{-})=H+2\frac{\p{x^i}}{\p{\tilde x^{-}}}A_i+
g\left(\frac{\p{x^i}}{\p{\tilde x^{-}}}\p_{i},\frac{\p{x^j}}{\p{\tilde x^{-}}}\p_{j}\right).$$ Hence,
$$q'=\tilde\p_{-}-\frac{1}{2}H'\p_{{+}}=q+U-\frac{1}{2}g(U,U)p,$$
where $$U=\frac{\p{x^i}}{\p{\tilde x^{-}}}X_i.$$
The equality $U=W$ is equivalent to  the system of equations
 \begin{equation}\label{eq1}\frac{\p{x^i}(\tilde x^{1},...,\tilde x^n,\tilde x^{-})}{\p{\tilde x^{-}}}=
 W^i(x^1(\tilde x^{1},..., \tilde x^n,\tilde x^{-}),...,x^n(\tilde x^{1},..., \tilde x^n,\tilde x^{-}),\tilde x^{-}).\end{equation}
Consider the system of ordinary differential equations
\begin{equation}\label{eq2}\frac{dy^i(\tilde x^{-})}{d \tilde x^{-}}=W^i(y^1(\tilde x^{-}),...,y^n(\tilde x^{-}),\tilde x^{-}).\end{equation}
Impose the initial conditions $y^i(\tilde x_0^{-})=\tilde x^{i}$.
Then for each set of numbers $\tilde x^{k}$ there exists a unique
solution $y^i(x^-)$. Since the solution depends smoothly on the
initial conditions, we may write the solution in the form
${x^i}(\tilde x^{1},...,\tilde x^n,\tilde x^{-})$. The obtained
functions satisfy  Equation \eqref{eq1}. Since $\det
\left(\frac{\partial x^i}{\partial \tilde x^{j}}(\tilde
x_0^{-})\right)\neq 0$, we get that $\det \left(\frac{\partial
x^i}{\partial \tilde x^{j}}\right)\neq 0$ for $\tilde x^{-}$ near
$\tilde x_0^{-}$. We obtain the required transformation. \eprf

We see that we may choose a Walker coordinate system such that
$v=0$ (if $v\neq 0$, take $W=-\frac{1}{\Lambda}v$, then $\tilde
v=0$). It can be shown that
$$v=-\left(\frac{1}{2}\p_{i}H_1-\Lambda A_i\right)h^{ij}X_j.$$
Since $\p_+\left(\left(\frac{1}{2}\p_{i}H_1-\Lambda
A_i\right)h^{ij}\right)=0$, it holds $\nabla_{\p_+}W=0$. Hence we
may find a coordinate system, where
$A_i=\frac{1}{2\Lambda}\p_{i}H_1$. Let us fix this system. In
\cite{G-P} it is noted that under the transformation
$$\tilde x^+=x^+-f(x^1,...,x^n,x^{-}),\quad \tilde x^i=x^i,\quad
\tilde x^{-}=x^{-}$$ the metric \eqref{Walker} changes in the
following way
 \begin{equation}\label{GP4.12} A_i\mapsto A_i+\partial_{i}f,\quad H_1\mapsto H_1+2\Lambda f,\quad H_0\mapsto H_0+H_1 f+\Lambda f^2+2\dot f.\end{equation}
Thus if we take $f=-\frac{1}{2\Lambda }H_1$, then with respect to
the new coordinates we have  $A_i=H_1=0$. Now  Theorem \ref{Main
Theorem} follows from (\ref{GP4.8}--\ref{GP4.11}). \hfill $\Box$

\subsection*{Proof of Proposition \ref{prop2}}
The decomposition of the $\lason$-projection of the holonomy as in
\eqref{LM0B}, $\h=\{0\}\oplus\h_1\oplus\cdots\oplus\h_r$ defines
parallel distributions $E^0, \ldots , E^r$, all containing the
parallel distribution of null lines. These distributions, in turn,
define coordinates
$$x^+,x_0^1,\ldots ,x_0^{n_0},\ldots ,x_r^1,...,x_r^{n_r},x^-$$  such that
$E^\alpha$ is spanned by $\p_+, \frac{\p}{\p x^1_\alpha},\ldots ,  \frac{\p}{\p x^{n_\alpha}_\alpha}$
and such that they are adapted in the sense of Section \ref{secConseq}.
Note that the most general coordinate transformation preserving these properties is given by
\begin{equation}
\begin{array}{rcl}
\tilde{x}^+&=&x^+ + \vf (x^1_0, \ldots , x^{n_r}_{r} ,x^-),\\
\tilde{x}^i_0 &=&\sum_{j=1}^{n_0}a^i_j x^j_0+b^i,\ \text{ for $
i=1 ,\ldots , n_0$,}
\\
\tilde{x}^i_\alpha &=&\psi^i_\alpha (x^1_\alpha, \ldots ,
x^{n_\alpha}_\alpha,x^-),\  \text{ for $ i=1 ,\ldots , n_\alpha$
and $\alpha=1, \ldots , r$,}
\\
\tilde{x}^-&=&x^-+c,
\end{array}\label{gaugeadapt}
\end{equation}
here $\tfrac{\p^2}{\p x^j_\beta \p x^i_\alpha}\vf=0$ if
$\beta\neq\alpha$, $(a^i_j)_{i,j=1}^{n_0}$ is an orthogonal matrix
and $b^i\in\Real$. Choosing $\vf\equiv 0$, it is clear that
Equation \eqref{ck1} can be solved separately for each $\alpha=1,
\ldots , r$. This shows that the coordinates found in Theorem
\ref{theo1} can be chosen to be adapted.

Now we turn to the second  statement of Proposition  \ref{prop2}.
Let us assume that $\Lambda\not=0$. Starting with adapted
coordinates, Equation \eqref{GP4.9} shows  that
\begin{equation} \label{mixedh1} \tfrac{\p^2}{\p x^j_\beta \p
x^i_\alpha}H_1=0,\ \ \text{ if $\beta\neq\alpha$.}\end{equation}
Consider the proof of Theorem \ref{Main Theorem} applied to a
metric in adapted coordinates in order to prove the second
statement. Equation \eqref{GP4.11} shows that $n_0=0$.   Recall
that we consider the system of equations \eqref{eq1} for
$W^i=\frac{1}{\Lambda}\left(\frac{1}{2}\p_jH_1-\Lambda
A_j\right)h^{ij}$. Since we have the property \eqref{mixedh1}, we
get that if the index $i$ corresponds to the space
$\Real^{n_\alpha}$, then $\frac{\p}{\p x^k_\beta} W^{i}=0$ if
$\beta\neq\alpha$. It is obvious that we get $r$ independent
systems of equations, each of these systems is a system  with
respect to the unknown
  functions $x^1_\alpha(\tilde x^1_\alpha,...,\tilde x^{n_\alpha}_\alpha),...,
x^{n_\alpha}_\alpha(\tilde x^1_\alpha,...,\tilde x^{n_\alpha}_\alpha)$.
It is clear that the  solution for such a system obtained above satisfies the requirements of the proposition. \hfill $\Box$

\subsection*{Proof of Proposition \ref{prop1}}
 As above, let $R=R(\lambda,v,R_0,P,T)$. Consider the coordinate system as in Theorem \ref{Main Theorem}. Then, $v=0$ and $\tRic(P)=0$.
The decomposition \eqref{LM0B} implies $P=P_1+\cdots + P_r$, where  $P_\beta\in\P(\h_\beta)$. Consequently, each $\tRic(P_\beta)$ is zero. Since
$\h_\alpha\subset\so(n_\alpha)$ is a symmetric Berger algebra, the equality $\tRic(P_\alpha)=0$ implies $P_\alpha=0$ \cite{GalP}, and this
is exactly Equation \eqref{4.9strong}. \hfill
$\Box$

\section{Examples}\label{secEx}

Suppose that metric \eqref{Walker} is Einstein with the cosmological constant $\Lambda\neq 0$. Then  \eqref{GP4.7} holds.
According to Theorem \ref{Main Theorem}, there exist new Walker coordinates $(\tilde x^+,\tilde x^1,...,\tilde x^n,\tilde x^-)$
such that $\tilde A=0$ and $\tilde H_1=0$. The proof of Theorem \ref{Main Theorem} implies that such coordinates can be found in the following way.
Consider the system of ordinary differential equations
\begin{equation}\label{eq2A}\frac{dy^i(\tilde x^{-})}{d \tilde x^{-}}=W^i(y^1(\tilde x^{-}),...,y^n(\tilde x^{-}),\tilde x^{-}),\end{equation} where $W^i=(\frac{1}{2\Lambda}\p_{j}H_1-
A_j)h^{ij}$ and impose the initial conditions $y^i(\tilde x_0^{-})=\tilde x^{i}$.
This will give the inverse transformation $$x^+=\tilde x^{+},\quad
x^i=x^i(\tilde x^{1},...,\tilde x^n,\tilde x^{-}),\quad  x^{-}=\tilde x^{-}$$
and allow to find the metric with respect to the new coordinates.
Note that $\tilde H_1=H_1$. If $H_1=0$, then with respect to the obtained coordinates $\tilde A_i=\tilde H_1=0$ holds.
If $H_1\neq 0$, then it is necessary to consider the additional transformation
$$\label{trans1}\tilde x^+\mapsto   \tilde x^++\frac{1}{2\Lambda}H_1,\quad \tilde x^i\mapsto \tilde x^i,\quad \tilde
x^{-}\mapsto \tilde x^{-}.$$
After this $\tilde A_i=\tilde H_1=0$.

The required coordinates can be found also in the following way.
First consider the transformation $$ x^+\mapsto   x^++\frac{1}{2\Lambda}H_1,\quad  x^i\mapsto  x^i,\quad
x^{-}\mapsto  x^{-}.$$ After this $H_1=0$ and $A_i$ changes to $A_i-\frac{1}{2\Lambda}\partial_iH_1$.
After this consider the system of ordinary differential equations \eqref{eq2A} with $W^i=-A_jh^{ij}$ and impose the initial conditions $y^i(\tilde x_0^{-})=\tilde x^{i}$.
With respect to the obtained coordinates $\tilde A_i=\tilde H_1=0$ holds.

\vskip0.3cm

For $n=2$ and $\Lambda\neq 0$ all solutions to Equation \eqref{Einstein} for metric \eqref{Walker} are obtained in \cite{Lewandowski}. It is proved that any such metric is given in the following way
 (we use slight modifications).
There exist coordinates $x^+,u,v,x^-$ such that $$g=\frac{2}{P^2}\d z\d\bar{z}+\left(2\d x^++2W\d z+2\bar{W}\d\bar{z}+\big(\Lambda\cdot (x^+)^2+H_0\big)\d x^-\right)\d x^-,$$
where \begin{align*} z&=u+iv,\quad 2P^2=|\Lambda|2P_0^2=|\Lambda|\left(1+\frac{\Lambda}{|\Lambda|} z\bar{z}\right)^2,\quad W=i\p_z L,\\
L&=2\text{\rm Re}\left(f\p_z(\ln P_0)-\frac{1}{2}\p_zf\right),\end{align*}
$f=f(z,x^-)$ is an arbitrary function holomorphic in $z$ and smooth in $x^-$, the function $H_0=H_0(z,\bar{z},x^-)$ can be expressed in a similar way in terms of $f$ and  another arbitrary function holomorphic in $z$ and smooth in $x^-$.

Using  this result, we consider several examples.

\begin{ex} Let $\Lambda<0$ and $f=c(x^-)$, we obtain the following metric
\begin{multline*}g=2\d x^+\d x^-+\frac{4}{-\Lambda\cdot(1-u^2-v^2)^2}\big((\d u)^2+(\d v)^2\big)\\+
\frac{c(x^-)}{(1-u^2-v^2)^2}\big(-4uv\d u +2(u^2-v^2+1)\d v\big)\d x^-+(\Lambda \cdot (x^+)^2+H_0)(\d x^-)^2,\end{multline*}
which becomes Einstein after a proper choice of the function $H_0$.
Equations \eqref{eq2A} take the form
$$\frac{\partial u}{\partial \tilde x^-}=-\frac{\Lambda}{2}uvc(x^-), \quad \frac{\partial v}{\partial \tilde x^-}=
\frac{\Lambda}{4}(u^2-v^2+1)c(x^-).$$
Using Maple 12, we find that the general solution of this system has the form
\begin{align*} u&=\frac{64c_1\Lambda^2}
{\left(c_1^2
\left(4 e^{-\frac{1}{2}\Lambda b(\tilde x^-)}+\Lambda c_2\right)^2+64\Lambda^4\right)e^{\frac{1}{2}\Lambda \tilde b(x^-)}},\\
v&=\frac{-16c_1^2e^{-\Lambda b(\tilde x^-)}+c_1^2c_2^2\Lambda^2
+64\Lambda^4} {c_1^2 \left(4 e^{-\frac{1}{2}\Lambda b(\tilde
x^-)}+\Lambda c_2\right)^2+64\Lambda^4},\end{align*} where $c_1$
and $c_2$ are arbitrary functions of $\tilde u$ and $\tilde v$,
$b(\tilde x^-)$ is the function such that $\frac{\d b(\tilde
x^-)}{\d \tilde x^-}=c(\tilde x^-) $ and $b(0)=0$. Substituting
the initial conditions $u(0)=\tilde u$,  $v(0)=\tilde v$, we
obtain $$c_1=\frac{\tilde u^2+\tilde v^2-2\tilde v^2+1}{\tilde
u}\Lambda^2,\quad c_2=-4\frac{\tilde u^2+\tilde
v^2-1}{\Lambda\cdot(\tilde u^2+\tilde v^2-2\tilde v^2+1)}.$$ With
respect to the obtained coordinates, we get
\begin{equation}g=2\d x^+\d x^-+\frac{4}{-\Lambda\cdot(1-u^2-v^2)^2}\big((\d u)^2+(\d v)^2\big)
+(\Lambda \cdot (x^+)^2+\tilde H_0)(\d x^-)^2.\end{equation} The
metric $g$ is Einstein if and only if $(\p_u^2+\p_v^2)\tilde
H_0=0$. Taking sufficiently general solutions of this equation
(e.g. $\tilde H_0=uv$), we obtain that this metric is
indecomposable and its  holonomy algebra is isomorphic to
$(\Real\oplus\so(2))\zr\Real^2$.

Note that taking $f=z^2$, one obtains the same example.
\end{ex}

\begin{ex} Let $\Lambda<0$ and $f=z c(x^-)$, we obtain the following metric
\begin{multline*}g=2\d x^+\d x^-+\frac{4}{-\Lambda\cdot(1-u^2-v^2)^2}\big((\d u)^2+(\d v)^2\big)\\+
\frac{2c(x^-)}{(1-u^2-v^2)^2}\big(v\d u-u\d v\big)\d x^-+(\Lambda \cdot (x^+)^2+H_0)(\d x^-)^2.\end{multline*}
Equations \eqref{eq2A} take the form
$$\frac{\partial u}{\partial \tilde x^-}=\frac{\Lambda}{4} vc(x^-), \quad \frac{\partial v}{\partial \tilde x^-}=
-\frac{\Lambda}{4} uc(x^-).$$
The general solution of this system has the form
$$u=c_1\cos\left(\frac{\Lambda}{4} b(\tilde x^-)\right)+c_2\sin\left(\frac{\Lambda}{4}b(\tilde x^-)\right),\quad v=-c_1\sin\left(\frac{\Lambda}{4}b(\tilde x^-)\right)+c_2\cos\left(\frac{\Lambda}{4}b(\tilde x^-)\right),$$
where $c_1$ and $c_2$ are arbitrary functions of $\tilde u$ and
$\tilde v$, and $b(\tilde x^-)$ is the function such that
$\frac{\d b(\tilde x^-)}{\d \tilde x^-}=c(\tilde x^-) $ and
$b(0)=0$. Substituting the initial conditions $u(0)=\tilde u$,
$v(0)=\tilde v$, we obtain $c_1=\tilde u$, $c_2=\tilde v$. With
respect to the obtained coordinates, we again get
\begin{equation}g=2\d x^+\d x^-+\frac{4}{-\Lambda\cdot(1-u^2-v^2)^2}\big((\d u)^2+(\d v)^2\big)+
(\Lambda \cdot (x^+)^2+\tilde H_0)(\d x^-)^2.\end{equation}
\end{ex}

\begin{ex} Let $\Lambda>0$ and $f=zc(x^-)$, we obtain the following metric
\begin{multline*}g=2\d x^+\d x^-+\frac{4}{\Lambda\cdot(1+u^2+v^2)^2}\big((\d u)^2+(\d v)^2\big)\\+
\frac{2c(x^-)}{(1+u^2+v^2)^2}\big(v\d u-u\d v\big)\d x^-+(\Lambda \cdot (x^+)^2+H_0)(\d x^-)^2.\end{multline*}
Equations \eqref{eq2A} take the form
$$\frac{\partial u}{\partial \tilde x^-}=-\frac{\Lambda}{4} vc(x^-), \quad \frac{\partial v}{\partial \tilde x^-}=
\frac{\Lambda}{4} uc(x^-).$$
The general solution of this system has the form
$$u=c_1\cos\left(\frac{\Lambda}{4}b(\tilde x^-)\right)+c_2\sin\left(\frac{\Lambda}{4}b(\tilde x^-)\right),\quad v=c_1\sin\left(\frac{\Lambda}{4}b(\tilde x^-)\right)-c_2\cos\left(\frac{\Lambda}{4}b(\tilde x^-)\right),$$
where $c_1$ and $c_2$ are arbitrary functions of $\tilde u$ and
$\tilde v$, and $b(\tilde x^-)$ is the function such that
$\frac{\d b(\tilde x^-)}{\d \tilde x^-}=c(\tilde x^-) $ and
$b(0)=0$. Substituting the initial conditions $u(0)=\tilde u$,
$v(0)=\tilde v$, we obtain $c_1=\tilde u$, $c_2=-\tilde v$. With
respect to the obtained coordinates, we get
\begin{equation}g=2\d x^+\d x^-+\frac{4}{\Lambda\cdot(1+u^2+v^2)^2}\big((\d u)^2+(\d v)^2\big)+
(\Lambda \cdot (x^+)^2+\tilde H_0)(\d x^-)^2.\end{equation} The metric $g$ is Einstein if and only if $(\p_u^2+\p_v^2)\tilde H_0=0$. Taking sufficiently general solution of this equation (e.g. $\tilde H_0=uv$), we obtain that this metric is indecomposable and its holonomy algebra is isomorphic to $(\Real\oplus\so(2))\zr\Real^2$.
\end{ex}

For most of the other functions $f$ Equations \eqref{eq1} and
their solutions become much more difficult. Further examples are
considered in \cite{GalEn2}. In particular, in \cite{GalEn2} there
are obtained examples such that the Riemannian part $h$ depends
non-trivially on the parameter $x^-$.

\vskip0.3cm

Consider the general Walker metric \eqref{Walker}. Theorem \ref{theo1} shows that there exist coordinates $(\tilde x^+,\tilde x^1,...,\tilde x^n,\tilde x^-)$ such that $\tilde A=0$. These coordinates can be found as in the proof of Theorem \ref{theo1} or in the following alternative way.

Consider the transformation given by the inverse one $x^+=\tilde x^{+}$,
$x^i=x^i(\tilde x^{1},...,\tilde x^n,\tilde x^{-}),$  $x^{-}=\tilde x^{-}$. It holds
$$\tilde\p_{+}=\p_{+}, \quad
\tilde\p_{i}=\frac{\p{x^j}}{\p{\tilde x^{i}}}\p_{j},\quad
\tilde\p_{-}=\frac{\p{x^i}}{\p{\tilde x^{-}}}\p_{i}+\p_{{-}}.$$
For the new Walker metric we get $$\tilde A_i=\frac{\p{x^j}}{\p{\tilde x^{i}}}\left(A_j+h_{jk}\frac{\p{x^k}}{\p{\tilde x^{-}}}\right).$$
Hence, if the equalities \begin{equation}\label{eq0}\frac{\p{x^i}}{\p{\tilde x^{-}}}=-A_jh^{ji}\end{equation} hold, then
$\tilde A_i=0.$ Impose the  conditions $x^i(\tilde x^1,...,\tilde x^n,\tilde x_0^{-})=\tilde x^{i}$.
Then for each set of numbers $\tilde x^{k}$ there exists a unique
solution $x^i(x^-)$ of the above system of equations. Since the solution depends smoothly on the
initial conditions, we may write the solution in the form
${x^i}(\tilde x^{1},...,\tilde x^n,\tilde x^{-})$. The obtained
functions satisfy  Equation \eqref{eq0}. Since $\det
\left(\frac{\partial x^i}{\partial \tilde x^{j}}(\tilde
x_0^{-})\right)\neq 0$, we get that $\det \left(\frac{\partial
x^i}{\partial \tilde x^{j}}\right)\neq 0$ for $\tilde x^{-}$ near
$\tilde x_0^{-}$. We obtain the required transformation.

\vskip0.2cm

Ricci-flat Walker metrics in dimension 4 are found in \cite{KG1,KG2}. They are of the form
\begin{equation}\label{metricKG} g=2\d x^+\d x^-+(\d u)^2+(\d v)^2+2 A_1 \d x\d x^-+(-(\p_uA_1)x^++H_0)(\d x^-)^2,\end{equation}
where $A_1$ and $H_0$ satisfy $\p_+A_1=\p_+ H_0=0$, \begin{align}\label{A10}
\p_u^2A_1+\p_v^2 A_1&=0,\\ \label{H00} \p_u^2H_0+\p_v^2
H_0&=2\p_{-}\p_u A_1-2A_1\p_u^2A_1-(\p_u A_1)^2+(\p_v
A_1)^2.\end{align}

Note that in order to get rid of the function $A_1$ it is enough to consider the
transformation with the inverse one $$x^+=\tilde x^+,\,\, u=f(\tilde u,\tilde v,\tilde x^-),\,\,v= \tilde v,\,\, x^-=\tilde x^-$$ such that the function $f$ satisfies the equation
\begin{equation}\label{eqf}\p_{-}f(\tilde u,\tilde v,\tilde x^-)=-A_1(f(\tilde u,\tilde v,\tilde x^-),\tilde v, \tilde x^-).\end{equation} Imposing the condition $f(\tilde u,\tilde v,0)=\tilde u$, we may consider the coordinates $\tilde u$  and $\tilde v$ as the parameters, then the obtained equation is an ordinary differential equation.

\begin{ex}
It is clear that  $A_1=uv$ and $H_0=\frac{1}{12}(u^4-v^4)$ are solutions of \eqref{A10} and \eqref{H00}.
We get the following Ricci-flat
metric:
\begin{equation} g=2\d x^+\d x^-+(\d u)^2+(\d v)^2+2uv\d u\d x^-
+\left(-vx^++ \frac{1}{12}(u^4-v^4)\right)(\d x^-)^2.\end{equation}
Equation \eqref{eqf} takes the form $$\p_-f(\tilde u,\tilde v,\tilde x^-)=-f(\tilde u,\tilde v,\tilde x^-)\tilde v$$ and it defines the transformation
$$\tilde x^+= x^+,\,\,\tilde u= u e^{vx^-} ,\,\,\tilde v= v,\,\,\tilde x^-= x^-.$$
  With respect to the obtained coordinates, we get
\begin{multline} g=2\d x^+\d x^-
+e^{-2vx^-}(\d u)^2-2ux^-e^{-2vx^-}\d u \d v+\left(1+u^2(x^-)^2e^{-2vx^-}\right)(\d v)^2
\\+\left(-vx^+-u^2v^2e^{-2vx^-}-\frac{1}{12}v^4+\frac{1}{12}u^4e^{-4x^-v}\right)(\d x^-)^2.
\end{multline}
The holonomy algebra of this metric equals to $(\Real\oplus\so(2))\zr\Real^2$.
\end{ex}

\begin{ex}

The functions  $A_1=e^u\cos v$ and $H_0=-\frac{1}{4}(1+2v\sin 2v)e^{2u}$ are solutions of \eqref{A10} and \eqref{H00}.
We get the following Ricci-flat
metric:
\begin{equation} g=2\d x^+\d x^-+(\d u)^2+(\d v)^2+2e^u\cos v\d u\d x^-
+\left(-x^+e^u\cos v-\frac{1}{4}(1+2v\sin 2v)e^{2u}\right)(\d x^-)^2.\end{equation}
Equation \eqref{eqf} takes the form $$\p_-f(\tilde u,\tilde v,\tilde x^-)=-e^{f(\tilde u,\tilde v,\tilde x^-)}\cos\tilde v$$ and it defines the transformation
$$\tilde x^+= x^+,\,\,\tilde u=-\ln\left(e^{-u}-x^-\cos v\right) ,\,\,\tilde v= v,\,\,\tilde x^-= x^-.$$
  With respect to the obtained coordinates, we get
\begin{multline} g=2\d x^+\d x^-+
\frac{1}{\left(x^-e^u\cos v+1\right)^2}\left((\d u)^2+2x^-e^{u}\sin v\d u \d v+\left(1+x^-e^{u}\right)(\d v)^2\right)
\\-\frac{1}{4\left(x^-e^u\cos v+1\right)^2}\left(4x^+\left(x^-\cos^2v+e^{-u}\cos v\right)+1+4\cos^2v+2v\sin2v\right)(\d x^-)^2.
\end{multline}
The holonomy algebra of this metric equals to $(\Real\oplus\so(2))\zr\Real^2$.
\end{ex}

\subsection*{Acknowledgments} We thank Helga Baum and D.\,V.\,Alekseevsky  for  discussions on
the topic of this paper. The first author was supported by the grant  201/09/P039 of the
Grant Agency of Czech Republic and by the grant MSM~0021622409 of the Czech Ministry of Education.


\centerline{\rule{8cm}{.5pt}}

\end{document}